\documentclass[11pt]{article}

\usepackage[margin=1in]{geometry}

\usepackage{amsmath}
\usepackage{amssymb}
\usepackage{amsthm}
\usepackage{mathrsfs}

\usepackage{natbib}

\usepackage{graphicx}
\usepackage{xcolor}
\usepackage{float}
\usepackage{tikz}
\usetikzlibrary{arrows.meta,positioning}
\usepackage{longtable}
\usepackage{booktabs}
\usepackage{multirow}
\usepackage{array}
\usepackage{placeins}

\usepackage{hyperref}

\title{Exact Output Steering with Vanishing Controls
in Weighted Shift Dynamics\\[0.5em]
\rule{0.85\textwidth}{0.5pt}}

\author{%
  Yunied Puig\\
  \texttt{puigdedios@gmail.com}
}

\makeatletter

\@ifundefined{theorem}{
  \newtheorem{theorem}{Theorem}[section]
}{}

\@ifundefined{lemma}{
  \newtheorem{lemma}[theorem]{Lemma}
}{}

\@ifundefined{proposition}{
  \newtheorem{proposition}[theorem]{Proposition}
}{}

\@ifundefined{corollary}{
  \newtheorem{corollary}[theorem]{Corollary}
}{}

\@ifundefined{definition}{
  \newtheorem{definition}[theorem]{Definition}
}{}

\@ifundefined{example}{
  
}{}

\@ifundefined{remark}{
  
}{}

\@ifundefined{conjecture}{
  
}{}

\@ifundefined{problem}{
  
}{}

\@ifundefined{question}{
  
}{}

\@ifundefined{examples}{
  
}{}

\makeatother

\begin{document}
\date{}
\maketitle

\begin{abstract}
Exploiting the rich dynamics of wrapped weighted shifts and their residualizations, we characterize exact projected steering with vanishing initial controls.
 For wrapped shifts, steering is characterized
by an expansive cycle multiplier and a residue separating the
source and control channels. For residual shifts, a distinct
spectral--toral mechanism emerges: steering through selected
conjugate planes is equivalent to radial expansion and an
arithmetic parity condition on phase resonances, yielding
bounded-gap steering times and quantitative control estimates. For structured machine-learning representations built from orbit dictionaries, we derive a two-scale gain condition for steering coefficients while preserving decoded outputs.
\end{abstract}

\section{Introduction}

Weighted backward shifts are highly structured model operators with rich dynamics. Beginning with the classical hypercyclic
examples of Rolewicz, they have played a central role in the
study of orbit density, recurrence, mixing, and other dynamical
properties of linear operators
\citep{Bayart2009}. Their significance lies
in the combination of simple algebraic structure and rich
iterative behavior: explicit formulas for their powers permit
precise dynamical classifications that are generally difficult
for arbitrary operators. However, the classical phenomena of
hypercyclicity and topological mixing rely essentially on
infinite-dimensional settings. No linear operator on a
finite-dimensional normed space can be hypercyclic.

This raises a natural question: which forms of dynamical
reachability remain available for finite-dimensional weighted
shifts? We investigate this question through projected steering.
Instead of requiring one uncontrolled orbit to visit arbitrary
states, we prescribe source and target outputs, allow a small
perturbation of the initial state through a fixed control
channel, and ask whether the target can be reached exactly
after autonomous evolution. The central issue is whether the
required perturbation can vanish as the steering horizon
increases. This formulation connects the orbit geometry of
linear dynamics with one-shot output controllability, while addressing a fundamentally different problem from classical finite-horizon reachability
with distributed inputs. Indeed, classical state-to-output and output-to-output controllability
have been studied through rank conditions and controllability
Gramians \citep{SarachikKreindler1965,DanhaneLoheacJungers2023}, 
while null controllability with vanishing control energy
has been investigated for time-distributed inputs \citep{pandolfi2013vanishing}.
Here, we instead restrict control to a single vanishing
perturbation of the initial state, require exact output
transitions along a common sequence of horizons tending
to infinity, and characterize the resulting steering
property through the cyclic and spectral--toral structure
of weighted shifts.

Let \(T\in\mathcal L(Z)\) act on a finite-dimensional real
or complex normed space \(Z\), and let \(Y\) be an output
space. We consider source and control embeddings
\(J_{\mathrm{out}},J_{\mathrm{ctl}}:Y\to Z\) and an output
map \(P:Z\to Y\) satisfying
\[
PJ_{\mathrm{out}}=I_Y,
\qquad
PJ_{\mathrm{ctl}}=0.
\]
Given arbitrary \(u,v\in Y\), we seek controls
\(h_k(u,v)\to0\) and a common increasing sequence
\(n_k\to\infty\) such that
\[
PT^{n_k}
\bigl(
J_{\mathrm{out}}u+
J_{\mathrm{ctl}}h_k(u,v)
\bigr)=v.
\]
We require linear dependence of the controls on the
source--target pair, with the corresponding operator
norms converging to zero. Thus exact output transitions
must be achievable uniformly on bounded source--target
sets while the initial intervention becomes negligible.

For a wrapped weighted shift, we obtain a complete characterization
for arbitrary fixed channels satisfying the stated geometry.
Steering is equivalent to cycle expansion, annihilation of the source
block at one residue, and invertibility of the corresponding control
block. An explicit channel construction attains the output dimension
allowed by the source--control geometry.

For real residual shifts with negative cycle multiplier,
we use the classical rotation--dilation decomposition and explicitly
constructed conjugate-plane channels. Steering for these channels
is equivalent to radial expansion and simultaneous angular recurrence. Kronecker's theorem \citep{Cassels1957} converts the recurrence
requirement into an exact arithmetic condition on the
resonance lattice. The resulting control formulas give amplitude bounds and bounded-gap
admissible horizons for every fixed positive control budget on bounded
source--target sets.

The resulting theory applies to structured machine-learning representations built from orbit dictionaries. Such representations arise in recent
orbit-based learning and structure-guided embedding
compression frameworks
\citep{PuigWBSNN2026,PuigOrBIT2026}.
Applying the steering criterion to a direct-sum lift, we derive a
sufficient expansion--contraction balance for preserving decoded
realizations across internal orbit depths.

\begin{definition}[The wrapped weighted backward shift]
Let
$Z=\mathbb F^{d}$,
$\
\mathbb F\in\{\mathbb R,\mathbb C\}$,
with canonical basis
$e_0,\ldots,e_{d-1}$.
A wrapped weighted backward shift $W:Z\to Z$ is defined as  
$We_j=w_je_{j-1\pmod{d}}$,
 with
$w_j\neq0$.
\end{definition}
\label{def:wrapped wbs}
Unlike the standard finite-dimensional backward shift, which is nilpotent, the wrapped shift allows indefinite iteration without collapsing to the zero operator, thereby preserving nontrivial dynamical behavior. The nonzero-weight assumption makes $W$ invertible.
Moreover, coordinates circulate modulo the ambient dimension instead of being discarded, and complete cycles produce scalar rescaling.
Define the full-cycle multiplier
$
\Omega:=\prod_{j=0}^{d-1}w_j.
$
Since every coordinate crosses every weight exactly
once during \(d\) iterations,
$
W^d=\Omega I.
$
Consequently, for \(n=qd+r\), with \(0\leq r<d\),
$
W^n=\Omega^qW^r.
$


\section{Projected Steering Geometry}

\paragraph{Output and control geometry.}
Let \(Z\) be a finite-dimensional normed space over
\(\mathbb F\in\{\mathbb R,\mathbb C\}\), let
\(T\in\mathcal L(Z)\), and set \(Y=\mathbb F^p\).
Fix injective linear maps
\(J_{\mathrm{out}},J_{\mathrm{ctl}}:Y\to Z\)
and an output map \(P:Z\to Y\) satisfying
\[
PJ_{\mathrm{out}}=I_Y,
\qquad
PJ_{\mathrm{ctl}}=0.
\]
Thus \(J_{\mathrm{out}}u\) represents the initial
source, while \(J_{\mathrm{ctl}}h\) is a perturbation
restricted to a fixed subspace of \(\ker P\).
For each \(n\geq0\), define the source-to-output
and control-to-output operators
\[
A_n:=PT^nJ_{\mathrm{out}},
\qquad
G_n:=PT^nJ_{\mathrm{ctl}}.
\]
Exact steering from \(u\in Y\) to \(v\in Y\)
at depth \(n\) amounts to solving
\[
A_nu+G_nh=v.
\]
Whenever \(G_n\) is invertible, the unique
control is
$
h=G_n^{-1}(v-A_nu).
$
The question is whether such controls can be
chosen to vanish along a common sequence of
increasing depths, uniformly for bounded
source--target pairs.

\begin{definition}[Steering Criterion]
\label{def:steering criterion}
We say that $T$ satisfies the Steering Criterion along a strictly increasing sequence $(n_k)$ if there exist linear maps $Q_k,S_k:Y\to Y$ such that
$
\|Q_k\|\to0,
\
\|S_k\|\to0,
$
and
\[
PT^{n_k}
\left(
J_{\mathrm{out}}+J_{\mathrm{ctl}}Q_k
\right)=0,
\qquad
PT^{n_k}J_{\mathrm{ctl}}S_k=I_Y.
\]
\end{definition}

Our steering criterion says, schematically,
\[
\forall u,v\quad
\exists (x_k(u,v))_k:
\quad
x_k(u,v)\to J_{\mathrm{out}}u,
\quad
PT^{n_k}x_k(u,v)=v,
\]
where $x_k(u, v)= J_{\mathrm{out}}u+J_{\mathrm{ctl}}(Q_{k}u+S_{k}v)$.

\begin{theorem}[Exact block characterization]
\label{thm:block charact}
Let $(n_k)$ be strictly increasing. The following are equivalent.

(i) $T$ satisfies the Steering Criterion along $(n_k)$;

(ii) For every $k$, $G_{n_k}$ is invertible and
$
\|G_{n_k}^{-1}\|\to0,
\
\|G_{n_k}^{-1}A_{n_k}\|\to0.
$
\end{theorem}

\begin{proof}
If (ii) holds, take
$
Q_k=-G_{n_k}^{-1}A_{n_k},
\
S_k=G_{n_k}^{-1}.
$
Then
$
A_{n_k}+G_{n_k}Q_k=0
$
and
$
G_{n_k}S_k=I_Y.
$
Conversely, if (i) holds, then
$
G_{n_k}S_k=I_Y.
$
Because $Y$ is finite-dimensional and $G_{n_k}$ is a square operator, $G_{n_k}$ is surjective and hence invertible. Therefore
$
S_k=G_{n_k}^{-1}.
$
Similarly,
$
A_{n_k}+G_{n_k}Q_k=0
$
implies
$
Q_k=-G_{n_k}^{-1}A_{n_k}.
$
The required limits are exactly $\|S_k\|\to0$ and $\|Q_k\|\to0$.
\end{proof}

\paragraph{Why projected rather than full-state steering?}
The Steering Criterion guarantees exact arrival at the prescribed output,
$PT^{n_k}x_k(u,v)=v$, whereas full-state asymptotic arrival would require
$T^{n_k}x_k(u,v)= J_{\mathrm{out}}v$.
We show in Appendix~\ref{app:full-state-obstruction} that, for both
expansive steering mechanisms classified below, every $u\neq0$
and every $v\in Y$ satisfy
\[
\|T^{n_k}x_k(u,v)-J_{\mathrm{out}}v\|
\longrightarrow\infty.
\]
In particular, neither steering mechanism can be asymptotically
full-state accurate.

\section{Cycle-Synchronous Steering dynamics for $W$}
\label{sec:steering for W}
Let \(W\) be a wrapped weighted backward shift, with full-cycle
multiplier \(\Omega\) as in Definition~\ref{def:wrapped wbs}. For the fixed channels,
write
\[
A_r=PW^rJ_{\mathrm{out}},
\qquad
G_r=PW^rJ_{\mathrm{ctl}},
\qquad 0\le r<d.
\]
The cycle identity gives
\[
A_{qd+r}=\Omega^qA_r,
\qquad
G_{qd+r}=\Omega^qG_r.
\]
Whenever \(G_r\) is invertible,
\[
G_{qd+r}^{-1}=\Omega^{-q}G_r^{-1},
\qquad
G_{qd+r}^{-1}A_{qd+r}=G_r^{-1}A_r.
\]
In particular, the normalized source block is independent of \(q\). For $R>0$, let $B_Y(0,R):=\{y\in Y:\|y\|<R\}$.

\begin{theorem}[Equivalent characterizations of steering for $W$]
\label{thm:characterization W}
Let $W$ be a wrapped weighted backward shift on $\mathbb F^d$,
and let
$
\Omega
:=
\prod_{j=0}^{d-1}w_j.
$
The following statements are equivalent.

\begin{enumerate}
\item[\rm(i)]
$W$ satisfies the Steering Criterion along some strictly increasing sequence.

\item[\rm(ii)]
One has
$
|\Omega|>1,
$
and there exists 
$
r\in\{0,\ldots,d-1\}
$
such that
$
A_r=0$ and $G_r$ is invertible.

\item[\rm(iii)]
There exists 
$
r\in\{0,\ldots,d-1\}
$
such that
$
A_r=0
$
and, for every $\varepsilon>0$ and every $R>0$, there exists an integer $q_0\geq0$ such that for every $q\geq q_0$,
\[
B_Y(0,R)
\subseteq
PW^{qd+r}J_{\mathrm{ctl}}
\bigl(B_Y(0,\varepsilon)\bigr).
\]
\end{enumerate}
\end{theorem}

\begin{proof}
Suppose (i) holds. By Theorem~\ref{thm:block charact}, along an increasing
sequence \((n_k)\), the blocks \(G_{n_k}\) are eventually
invertible and
\[
\|G_{n_k}^{-1}\|\to0,
\qquad
\|G_{n_k}^{-1}A_{n_k}\|\to0.
\]
Passing to an infinite subsequence with a fixed residue
\(r\), write \(n_k=q_kd+r\), where \(q_k\to\infty\).
The preceding identities imply that \(G_r\) is
invertible, \(|\Omega|>1\), and
\(G_r^{-1}A_r=0\). Thus \(A_r=0\), proving (ii).
Conversely, assume (ii) and take \(n_k=kd+r\).
Then
\[
A_{n_k}=0,
\qquad
G_{n_k}^{-1}=\Omega^{-k}G_r^{-1}.
\]
Hence the Steering Criterion holds by Theorem~\ref{thm:block charact},
with explicit maps
\[
Q_k=0,
\qquad
S_k=\Omega^{-k}G_r^{-1}.
\]

Finally, if (ii) holds, then
\(PW^{qd+r}J_{\mathrm{ctl}}=\Omega^qG_r\).
For every \(\varepsilon>0\),
\[
B_Y\left(
0,\frac{\varepsilon|\Omega|^q}{\|G_r^{-1}\|}
\right)
\subseteq
PW^{qd+r}J_{\mathrm{ctl}}B_Y(0,\varepsilon).
\]
Since \(|\Omega|>1\), this proves (iii).
Conversely, (iii) forces \(G_r\) to be surjective,
hence invertible, because its range contains a
nonempty open ball. If \(|\Omega|\leq1\), then
\[
PW^{qd+r}J_{\mathrm{ctl}}B_Y(0,1)
\subseteq B_Y(0,\|G_r\|)
\]
for every \(q\), contradicting (iii) for sufficiently
large target balls. Thus (ii) follows.
\end{proof}

\begin{lemma}[Maximal steerable output dimension for $W$]
\label{lemm:max steerable dim for W}
If \(Z=\mathbb F^d\), \(Y=\mathbb F^p\) and $W$ satisfies the Steering Criterion along some sequence, then  $p\le\left\lfloor\frac d2\right\rfloor.$
\end{lemma}
\begin{proof}
By Theorem \ref{thm:characterization W} there exists $0\le r\le d-1$ such that $A_r=PW^rJ_{\mathrm{out}}=0$, hence $W^rJ_{\mathrm{out}}Y\subseteq\ker P$.
But $W^r$ and $J_{\mathrm{out}}$ are injective, so
$
\dim W^rJ_{\mathrm{out}}Y=p.
$
Since
$
\dim\ker P=d-p,
$
we obtain
$
p\le d-p.
$
Therefore
$
p\le\left\lfloor\frac d2\right\rfloor.
$
\end{proof}

\paragraph{Canonical steering channels for the wrapped shift.}
The output and control channels can be constructed explicitly. 

\begin{proposition}[Canonical steering geometry]
\label{prop:canonical steering geometry}
Let
$W :\mathbb F^{d}\to
\mathbb F^{d}$
be a wrapped weighted backward shift with nonzero weights, and assume
$d\geq2p$.
Let $f_0,\ldots,f_{p-1}$ be the canonical basis of $Y=\mathbb F^p$, and let $e_0,\ldots,e_{d-1}$ be the canonical basis of $Z$.
Define
$J_{\mathrm{out}}f_j:=e_j$,
$J_{\mathrm{ctl}}f_j:=e_{p+j}$,
and
\[
Pe_j=
\begin{cases}
f_j, & 0\leq j\leq p-1,\\
0, & p\leq j\leq d-1.
\end{cases}
\]
Then
$PJ_{\mathrm{out}}=I_Y$ and $PJ_{\mathrm{ctl}}=0$.
Moreover, for the residue $r=p$,
$A_p=PW^pJ_{\mathrm{out}}=0$
and
$G_p=PW^pJ_{\mathrm{ctl}}$
is invertible.
\end{proposition}

\begin{proof}
For \(0\leq j<p\), the vector
\(W^pJ_{\mathrm{out}}f_j\)
is supported in the last \(p\) coordinates,
so \(PW^pJ_{\mathrm{out}}=0\).
On the other hand,
\(W^pJ_{\mathrm{ctl}}f_j\)
is a nonzero multiple of \(e_j\).
Thus \(G_p=PW^pJ_{\mathrm{ctl}}\)
is diagonal with nonzero diagonal entries, hence invertible.
\end{proof}

Consequently, when \(|\Omega|>1\), every output dimension
\(1\le p\le\lfloor d/2\rfloor\) is attained by channels satisfying
the Steering Criterion.

\section{Spectral--Toral Steering for \texorpdfstring{$R=I+W$}{R = I + W}}
\label{sec:steering for R}

Let \(W:\mathbb R^d\to\mathbb R^d\) be a real
wrapped weighted backward shift with
$
\Omega=-r^d,\ r>0,
$
and set \(R=I+W\). Write
\[
M=\left\lfloor\frac d2\right\rfloor,
\qquad
\delta=d-2M\in \{0, 1\}.
\]
When $d$ is even, $\delta=0$ and the whole real space decomposes into $M$ invariant two-dimensional conjugate-eigenvalue planes. When $d$ is odd, $\delta=1$ and there are $M$ such planes together with one additional real eigendirection.

The characteristic polynomial of \(W\) is
\(z^d+r^d\). Its eigenvalues are
\[
\nu_j=r\exp\left(\frac{(2j+1)\pi i}{d}\right),
\qquad 0\leq j<d.
\]
For \(0\leq j<M\), define
\[
\theta_j=\frac{(2j+1)\pi}{d},
\qquad
\lambda_j=1+re^{i\theta_j}
=\rho_je^{i\phi_j},
\]
where
$
\rho_j=\sqrt{1+r^2+2r\cos\theta_j}.
$
The nonreal eigenvalues occur in conjugate pairs,
and there is a real isomorphism \(\mathcal U\) such that
\[
\mathcal U^{-1}R^n \mathcal U=
\left(
\bigoplus_{j=0}^{M-1}
\rho_j^n\operatorname{Rot}(-n\phi_j)
\right)
\oplus[(1-r)^n]^\delta.
\]
The characteristic-polynomial calculation,
the similarity, and the construction of \(\mathcal U\)
are given in Appendix~\ref{app:spectral toral for R}.

Choose a nonempty set of conjugate-plane indices
\[
J=\{j_1,\ldots,j_p\}
\subseteq\{0,\ldots,M-1\}.
\]

Write
\[
\widehat Z=
\left(\bigoplus_{j=0}^{M-1}\mathbb R^2\right)
\oplus\mathbb R^\delta.
\]
For each selected plane \(j_\ell\), let
\(\widehat J_{\mathrm{out}}f_\ell\) and
\(\widehat J_{\mathrm{ctl}}f_\ell\) be the vectors
supported on that plane with respective coordinates
\((1,0)\) and \((0,1)\). Define
\[
\widehat P((a_j,b_j)_{j=0}^{M-1},t)
=
(a_{j_1},\ldots,a_{j_p}),
\]
where \(t\) is omitted when \(d\) is even.
Transporting these maps through the spectral
similarity \(U\), set
\[
J_{\mathrm{out}}=\mathcal U\widehat J_{\mathrm{out}},
\qquad
J_{\mathrm{ctl}}=\mathcal U\widehat J_{\mathrm{ctl}},
\qquad
P=\widehat P \mathcal U^{-1}.
\]
Then
\[
PJ_{\mathrm{out}}=I_Y,
\qquad
PJ_{\mathrm{ctl}}=0,
\qquad Y=\mathbb R^p.
\]
All details in Appendix \ref{app:explicit construction channels for R}. Equip \(Y=\mathbb R^p\) with its Euclidean norm.
Now, direct calculation (see Lemma \ref{lem:rotational source control blocks}) gives
\[
A_n=
\operatorname{diag}_{1\leq\ell\leq p}
\bigl(\rho_{j_\ell}^{n}
\cos(n\phi_{j_\ell})\bigr),
\quad
\text{ and }
\quad
G_n=
\operatorname{diag}_{1\leq\ell\leq p}
\bigl(\rho_{j_\ell}^{n}
\sin(n\phi_{j_\ell})\bigr).
\]
Whenever \(G_n\) is invertible,
\begin{equation}
\label{eq:translation of inverses toral case}
\|G_n^{-1}A_n\|
=
\max_{1\leq\ell\leq p}
|\cot(n\phi_{j_\ell})|,
\quad
\text{ and }
\quad
\|G_n^{-1}\|
=
\max_{1\leq\ell\leq p}
\rho_{j_\ell}^{-n}
|\csc(n\phi_{j_\ell})|.
\end{equation}
Thus the Steering Criterion separates into
simultaneous angular recurrence and radial expansion.
Set
\[
\alpha_J=
\left(
\frac{\phi_{j_1}}{\pi},
\ldots,
\frac{\phi_{j_p}}{\pi}
\right)\in\mathbb T^p,
\qquad
\beta=\left(\frac12,\ldots,\frac12\right),
\]
and define
$
H_{\alpha_J}
=
\overline{\{n\alpha_J:n\in\mathbb Z\}}
\subseteq\mathbb T^p.
$
 These identities lead to the following
exact classification.
Remember that a set $\mathcal N \subseteq \mathbb N$ is called syndetic if there exists an integer
$
L\geq1
$
such that for every $N\ge 1$, every interval
$
\{N,N+1,\ldots,N+L\}\subset \mathbb N,
$
contains at least one element of $\mathcal N$.

\begin{theorem}[Exact spectral--toral classification for $R$]
\label{thm:exact spectral toral classification R}
Let
$
R=I+W
$
on
$
Z=\mathbb R^d,
$
where $W$ is a real wrapped weighted backward shift with nonzero weights and full-cycle multiplier
$
\Omega=-r^d,
\
r>0.
$
Let
$
M=\left\lfloor\frac d2\right\rfloor,
$
choose
$
J=\{j_1,\ldots,j_p\}
\subseteq
\{0,\ldots,M-1\},
$
and define the channels
$
P,
\
J_{\mathrm{out}},
$ and $
J_{\mathrm{ctl}}
$
through the selected conjugate planes as above.

For every $\ell=1,\ldots,p$, define
$
\theta_{j_\ell}
=
\frac{(2j_\ell+1)\pi}{d},
\
\lambda_{j_\ell}
=
1+re^{i\theta_{j_\ell}}
=
\rho_{j_\ell}e^{i\phi_{j_\ell}},
$
and
$
\alpha_J
=
\left(
\frac{\phi_{j_1}}{\pi},
\ldots,
\frac{\phi_{j_p}}{\pi}
\right)
\in\mathbb T^p.
$
Then the following are equivalent:

\begin{enumerate}
\item[\rm(i)]
$R$ satisfies the Steering Criterion along some strictly increasing sequence
$
n_k\longrightarrow\infty
$
with respect to the selected conjugate-plane channels.

\item[\rm(ii)]
The selected radial factors 
$
\rho_{j_\ell}>1
\
\forall\ell=1,\ldots,p,
$
and
$
\beta
=
\left(
\frac12,\ldots,\frac12
\right)
\in H_{\alpha_J}.
$

\item[\rm(iii)]
The selected radial factors 
$
\rho_{j_\ell}>1
\
\forall\ell=1,\ldots,p,
$
and the parity-resonance condition holds:
$q\cdot\alpha_J\in\mathbb Z
\quad\Longrightarrow\quad
\sum_{\ell=1}^p q_\ell
\text{ is even},
\
\forall q\in\mathbb Z^p.
$
\item[\rm(iv)]
The selected radial factors 
$
\rho_{j_\ell}>1
\
\forall \ell=1,\ldots,p,
$
and, for every $\varepsilon>0$, the phase-recurrence set
$
\mathcal N_\varepsilon
:=
\left\{
n\geq1:
\|n\alpha_J-\beta\|_{\mathbb T^p}<\varepsilon
\right\}
$
is syndetic in $\mathbb N$.
\end{enumerate}
\end{theorem}

\textbf{Proof sketch.}
By Theorem~\ref{thm:block charact}, steering along \((n_k)\) is
equivalent to
\[
\|G_{n_k}^{-1}\|\to0,
\qquad
\|G_{n_k}^{-1}A_{n_k}\|\to0.
\]
Equation \ref{eq:translation of inverses toral case} shows that the second
limit holds precisely when $\cot(n_k\phi_{j_\ell})\to0
\ \text{for every }\ell$, which is equivalent to
$
n_k\alpha_J\to\beta \
\ \text{in }\mathbb T^p
$ by Lemma~\ref{lem:phase convergence cotangent}.
Such an increasing sequence exists if and only if
\(\beta\in H_{\alpha_J}\).
Along this sequence,
$
|\sin(n_k\phi_{j_\ell})|\to1
$
for every selected mode. Hence the first block
limit holds if and only if
$
\rho_{j_\ell}>1,
\ 1\leq\ell\leq p.
$
This proves the equivalence of (i) and (ii).
Kronecker's closure theorem \citep{Cassels1957}, gives
\[
H_{\alpha_J}
=
\left\{
x\in\mathbb T^p:
q\cdot x\in\mathbb Z
\text{ whenever }
q\cdot\alpha_J\in\mathbb Z
\right\}.
\]
Since
$
q\cdot\beta=\frac12\sum_{\ell=1}^p q_\ell,
$
conditions (ii) and (iii) are equivalent.
Finally, translation by \(\alpha_J\) is minimal
on the compact group \(H_{\alpha_J}\).
Therefore every nonempty relatively open
neighborhood of \(\beta\) has syndetic return
times. Conversely, syndetic returns to
arbitrarily small neighborhoods imply
\(\beta\in H_{\alpha_J}\).
This proves the equivalence with (iv).

This establishes the equivalence of the four conditions.
The complete proof, including the auxiliary spectral
and recurrence arguments, is given in
Appendix~\ref{app:detailed proof of exact spectral toral classification for R}.
\hfill$\blacksquare$

\paragraph{Expansion and reachability do not suffice.}
Consider $r=\sqrt{3}$ and
\[
W_r=
\begin{pmatrix}
0&r\\
-r&0
\end{pmatrix},
\qquad R_r=I+W_r,
\]
with $J_{\mathrm{out}}u=(u,0)$,
$J_{\mathrm{ctl}}h=(0,h)$, and $P(x_0,x_1)=x_0$.
Here $\rho=\sqrt{1+r^2}=2>1$.
Since $PR_r(J_{\mathrm{out}}u+J_{\mathrm{ctl}}h)
=u+rh$, every target $v$ is reachable at horizon one
by choosing $h=(v-u)/r$.
Nevertheless, $\phi=\arctan r=\pi/3$, and whenever
$G_n$ is invertible,
\[
\|G_n^{-1}A_n\|=|\cot(n\phi)|=\frac{1}{\sqrt3}.
\]
Thus the source-cancellation controls cannot vanish,
and the Steering Criterion fails despite radial expansion
and finite-time output reachability.
Equivalently, $\alpha=\phi/\pi=1/3$ violates the
parity-resonance condition of Theorem~\ref{thm:exact spectral toral classification R}, since
$3\alpha=1\in\mathbb Z$ while $3$ is odd.
This two-dimensional counterexample shows that the
arithmetic condition in Theorem~\ref{thm:exact spectral toral classification R} captures a genuine
obstruction not detected by expansion and reachability alone.

\paragraph{Control budgets and admissible horizons.}
For any fixed positive control budget and bounded source–target pairs, exact steering is possible at horizons with bounded gaps. Assume the equivalent conditions of Theorem~\ref{thm:exact spectral toral classification R}, and set
\(\rho_{\min}=\min_{1\le\ell\le p}\rho_{j_\ell}>1\).
For \(0<\varepsilon<1/2\) and \(n\in\mathcal N_\varepsilon\),
the unique exact steering control satisfies
\[
\|h_n(u,v)\|
\le
\tan(\pi\varepsilon)\|u\|
+
\frac{\rho_{\min}^{-n}}{\cos(\pi\varepsilon)}\|v\|.
\]
Indeed, phase proximity gives
\[
|\cot(n\phi_{j_\ell})|\le\tan(\pi\varepsilon),
\qquad
|\csc(n\phi_{j_\ell})|
\leq\frac{1}{\cos(\pi\varepsilon)},
\]
and the estimate follows from the diagonal control formulas (\ref{eq:inverse diag control formulas for R}).

Given \(B>0\) and \(\eta>0\), choose \(\varepsilon\) so that
\(B\tan(\pi\varepsilon)<\eta/2\).
For all sufficiently large \(n\in\mathcal N_\varepsilon\), the
second term is also below \(\eta/2\) whenever \(\|v\|\le B\).
Thus the horizons admitting exact controls of norm at most \(\eta\)
uniformly for \(\|u\|,\|v\|\le B\) form a syndetic set.
No uniform gap bound as \(\eta\downarrow0\) is asserted.

The radial expansivity condition in Theorem~\ref{thm:exact spectral toral classification R}
admits an explicit threshold when all conjugate
rotation--dilation planes are selected. Set
\[
\kappa_d=
\begin{cases}
1, & d \text{ is even},\\
2, & d \text{ is odd}.
\end{cases}
\]
Then
$
\min_{0\leq j<M}\rho_j^2
=
1+r^2-2r\cos\left(\frac{\kappa_d\pi}{d}\right).
$
Consequently, all conjugate residual modes are
expansive if and only if
$
r>2\cos\left(\frac{\kappa_d\pi}{d}\right).
$
The proof is given in Appendix~\ref{app:smallest residual radial factor}.

\section{Steering Orbit-Dictionary Representations}
\label{sec:framework-navigation}
Orbit dictionaries yield structured machine-learning representations using operator-generated atoms. We apply our steering criteria to their coefficients and derive a gain condition preserving decoded realizations across depths.
The construction applies,
in particular, to the orbit-based representations considered in
\citep{PuigWBSNN2026} and \citep{PuigOrBIT2026}.


\paragraph{Finite orbit representations.}
Let \(X=\mathbb F^{d_{\mathrm{scaff}}}\),
\(V=\mathbb F^C\), and let
\(W_{\mathrm{scaff}}\in\mathcal L(X)\) be a
wrapped weighted shift. Fix finite blocks \(D_1,\ldots,D_K\), vectors
\(x_{k,j}\in X\), \(1\le j\le |D_k|\), and linear maps
\(J_k:X\to V\).
Set
\[
N_{\mathrm{sub}}=\sum_{k=1}^K|D_k|,
\qquad
Y_{\mathrm{fr}}
=\mathbb F^{d_{\mathrm{scaff}}N_{\mathrm{sub}}}.
\]
For \(c=(c_{k,m,j})\in Y_{\mathrm{fr}}\), define
the depth-dependent synthesis maps
\[
\mathscr D_s(c)=
\sum_{k=1}^K\sum_{j=1}^{|D_k|}
\sum_{m=0}^{d_{\mathrm{scaff}}-1}
c_{k,m,j}J_kW_{\mathrm{scaff}}^{m+s}x_{k,j}.
\]
In particular, \(\mathscr D_0\) is the synthesis map
of the fixed orbit dictionary.

\paragraph{Commuting steering dynamics.}
Let \(T_{\mathrm{steer}}\in\mathcal L(Z)\)
satisfy the Steering Criterion on
\(Y_{\mathrm{fr}}\). On
\(\mathcal Z=Z\oplus X^{N_{\mathrm{sub}}}\),
define
\[
\mathbb T
=T_{\mathrm{steer}}\oplus I,
\qquad
\mathbb W
=I_Z\oplus
W_{\mathrm{scaff}}^{\oplus N_{\mathrm{sub}}}.
\]
These operators commute by construction, since they act on separate
summands. Writing
\(\mathbf x=(x_{k,j})_{k,j}\), define
the decoder
\[
\mathfrak D(z,\mathbf x)=
\sum_{k,j,m}
(Pz)_{k,m,j}J_kW_{\mathrm{scaff}}^m x_{k,j},
\]
where the indices range over the finite
orbit dictionary. 
Theorem~\ref{thm:characterization W} with Proposition~\ref{prop:canonical steering geometry}, or Theorem~\ref{thm:exact spectral toral classification R}, supplies explicit steering
operators under the corresponding hypotheses.
By definition, for every \(z\in Z\) and \(s\ge0\),
\[
\mathfrak D\!\left(
z,(W_{\mathrm{scaff}}^{\oplus N_{\mathrm{sub}}})^s\mathbf x
\right)
=
\mathscr D_s(Pz).
\]

\begin{theorem}[Controlled framework navigation]
\label{thm:framework-navigation}
Suppose \(T_{\mathrm{steer}}\) satisfies
the Steering Criterion along \(n_q\to\infty\),
with steering maps \(Q_q,S_q\).
Fix \(c_0\in Y_{\mathrm{fr}}\), and let
\((c_q)\subset Y_{\mathrm{fr}}\) satisfy $\|S_q c_q\|\to 0$. For arbitrary internal
depths \(s_q\geq0\), set
\[
h_q=Q_qc_0+S_qc_q,
\qquad
\xi_q=
(J_{\mathrm{out}}c_0+J_{\mathrm{ctl}}h_q,
\mathbf x).
\]
Then \(h_q\to0\),
\(\xi_q\to(J_{\mathrm{out}}c_0,\mathbf x)\),
and
\begin{equation}
\label{eq:control navigation}
\mathfrak D(\mathbb T^{n_q}\mathbb W^{s_q}\xi_q)
=
\mathfrak D(\mathbb W^{s_q}\mathbb T^{n_q}\xi_q)
=
\mathscr D_{s_q}(c_q).
\end{equation}
\end{theorem}

\begin{proof}
Since the Steering Criterion guarantees \( \ \|Q_q\|\to0\), our assumption $\|S_q c_q\|\to 0$ implies
$
\|h_q\|
\leq
\|Q_q\|\|c_0\|+\|S_qc_q\|
\longrightarrow0.
$
Since $J_{\mathrm{ctl}}$ is bounded,
$J_{\mathrm{ctl}}h_q\to0$, and therefore
$\xi_q\to(J_{\mathrm{out}}c_0,\mathbf x)$. The Steering Criterion gives,
$
PT_{\mathrm{steer}}^{n_q}
(J_{\mathrm{out}}c_0+J_{\mathrm{ctl}}h_q)
=
PT_{\mathrm{steer}}^{n_q}
(J_{\mathrm{out}}+J_{\mathrm{ctl}}Q_q)c_0
+
PT_{\mathrm{steer}}^{n_q}
J_{\mathrm{ctl}}S_qc_q
=
0+c_q=c_q.
$
The commutation of the two lifted operators gives
\[
\begin{aligned}
\mathfrak D
\left(
\mathbb T_{\mathrm{steer}}^{n_q}
\mathbb W_{\mathrm{scaff}}^{s_q}\xi_q
\right)
&=
\sum_{k=1}^K
\sum_{m=0}^{d_{\mathrm{scaff}}-1}
\sum_{j=1}^{|D_k|}
\left[
PT_{\mathrm{steer}}^{n_q}
(J_{\mathrm{out}}c_0+J_{\mathrm{ctl}}h_q)
\right]_{k,m,j}
J_kW_{\mathrm{scaff}}^{m+s_q}x_{k,j}\\
&\quad=
\sum_{k=1}^K
\sum_{m=0}^{d_{\mathrm{scaff}}-1}
\sum_{j=1}^{|D_k|}
(c_q)_{k,m,j}
J_kW_{\mathrm{scaff}}^{m+s_q}x_{k,j}
\quad=
\mathscr D_{s_q}(c_q).
\end{aligned}
\]
\end{proof}

Every bounded target-code sequence satisfies
\(\|S_qc_q\|\to0\).

\paragraph{Preserving realizations across depths.}
Let \(\Omega_{\mathrm{scaff}}\neq0\)
be the full-cycle multiplier of
\(W_{\mathrm{scaff}}\).
Define the coefficient shift
\(C_{\mathrm{scaff}}\) by
\[
(C_{\mathrm{scaff}}c)_{k,m,j}
=
\begin{cases}
c_{k,m+1,j},
&0\leq m<d_{\mathrm{scaff}}-1,\\
\Omega_{\mathrm{scaff}}^{-1}c_{k,0,j},
&m=d_{\mathrm{scaff}}-1.
\end{cases}
\]
The wrapped-cycle identity implies
\begin{equation}
\label{eq:exact preservation of the same realization across depths}
\mathscr D_s(C_{\mathrm{scaff}}^sc)=\mathscr D_0(c),
\qquad
C_{\mathrm{scaff}}^{d_{\mathrm{scaff}}}
=\Omega_{\mathrm{scaff}}^{-1}I.
\end{equation}

 Proof of (\ref{eq:exact preservation of the same realization across depths}) in Appendix \ref{app:Contragredient Coefficient Transport}. Identities in (\ref{eq:exact preservation of the same realization across depths}) show that
internal orbit evolution can be compensated exactly,
preserving the same depth-zero realization across
arbitrary scaffold depths. However, when
$0<|\Omega_{\mathrm{scaff}}|<1$, this preservation
requires coefficient amplification as the depth
increases, potentially preventing the steering
controls from vanishing. The following result
establishes an explicit tradeoff between internal
scaffold depth and external steering horizon,
guaranteeing exact preservation of the realization
with asymptotically vanishing initial controls.

\begin{corollary}[Two-scale gain condition]
\label{cor:two-scale-gain}
In the setting of Theorem~\ref{thm:framework-navigation},
fix \(c_*\in Y_{\mathrm{fr}}\), and let
$
s_q=a_qd_{\mathrm{scaff}}+t_q,
\
0\leq t_q<d_{\mathrm{scaff}}.
$
Choose \(c_0=c_*\) and
\(c_q=C_{\mathrm{scaff}}^{s_q}c_*\).
Consequently, the fixed realization
\(\mathscr D_0(c_*)\) is preserved exactly with
vanishing initial controls whenever
\begin{equation}
\label{eq:abstract two scale gain condition}
\|S_q\||\Omega_{\mathrm{scaff}}|^{-a_q}
\longrightarrow0.
\end{equation}
For cycle-synchronous steering with \(n_q=q\,d_{\mathrm{steer}}+r\), it is
sufficient that
$
|\Omega_{\mathrm{steer}}|^q
|\Omega_{\mathrm{scaff}}|^{a_q}
\longrightarrow\infty,
$
whereas spectral--toral steering admits
the sufficient condition
$
\rho_{\min}^{n_q}
|\Omega_{\mathrm{scaff}}|^{a_q}
\longrightarrow\infty
$
along the recurrent steering horizons.
\end{corollary}

\begin{proof}
By the contragredient identities (\ref{eq:exact preservation of the same realization across depths}),
\[
\mathscr D_{s_q}(c_q)=\mathscr D_0(c_*),
\qquad
c_q=\Omega_{\mathrm{scaff}}^{-a_q}
C_{\mathrm{scaff}}^{t_q}c_*.
\]
Hence
\[
\|S_q c_q\|
\leq
M_{\mathrm{scaff}}\|S_q\|
|\Omega_{\mathrm{scaff}}|^{-a_q}\|c_*\|,
\]
where
$
M_{\mathrm{scaff}}
=
\max_{0\leq t<d_{\mathrm{scaff}}}
\|C_{\mathrm{scaff}}^t\|.
$ Theorem~\ref{thm:framework-navigation} then applies whenever (\ref{eq:abstract two scale gain condition}) is satisfied, yielding
\[
\mathfrak D
\left(
\mathbb T_{\mathrm{steer}}^{n_q}
\mathbb W_{\mathrm{scaff}}^{s_q}\xi_q
\right)
=
\mathscr D_{s_q}(c_q)
=
\mathscr D_0(c_\ast).
\]

To conclude, remember that for cycle-synchronous steering 
$\|S_q\|
=
|\Omega_{\mathrm{steer}}|^{-q}\|G_r^{-1}\|$,
while for the spectral--toral steering case, for all sufficiently large $q$,
$\|S_q\|\leq2\rho_{\min}^{-n_q}$.
\end{proof}


Extending this construction to redundant orbit charts, including
OrBIT, requires compatibility with chart gluing and is left for
future work.

\section{Conclusion}
We have developed a constructive mathematical characterization of exact output steering with vanishing initial controls for weighted shifts and selected residual-shift dynamics, connecting operator dynamics, arithmetic recurrence, control theory, and structured machine-learning representations.


\bibliographystyle{plainnat}
\bibliography{references}

\newpage

\section*{Appendix}

\appendix

\section{Spectral-Toral Steering Dynamics for \texorpdfstring{$R=I+W$}{R = I + W}}
\label{app:spectral toral for R}
The following elementary spectral calculations and classical
recurrence facts are included to make the channel construction
and Theorem~\ref{thm:exact spectral toral classification R} self-contained.

\begin{lemma}[Characteristic polynomial and eigenvalues of $W$]
\label{lem:spectrum of W negative cycle}
The minimal and characteristic polynomials of $W$ are both equal to
$
\mu_W(z)
=
\chi_W(z)
=
z^d+r^d.
$
Consequently, the eigenvalues of $W$ are simple and equal to
\[
\nu_\ell
=
r\exp\left(
\frac{(2\ell+1)\pi i}{d}
\right),
\qquad
\ell=0,\ldots,d-1.
\]
\end{lemma}

\begin{proof}
Since $W^d=-r^dI_{\mathbb R^d}$, the minimal
polynomial $\mu_W$ divides $z^d+r^d$.
On the other hand, the vectors
\[
e_0,We_0,\ldots,W^{d-1}e_0
\]
are nonzero multiples of distinct canonical basis
vectors, since each iteration shifts the coordinate
index by one modulo $d$ and all weights are nonzero.
They are therefore linearly independent, so no
nonzero polynomial of degree less than $d$
annihilates $W$. Thus $\deg\mu_W=d$, and hence
\[
\mu_W(z)=z^d+r^d.
\]
Since the characteristic polynomial $\chi_W$
is monic of degree $d$ and divisible by $\mu_W$,
we obtain $\chi_W=\mu_W$.
The eigenvalues are consequently the $d$ distinct
roots of $z^d=-r^d$, namely
\[
\nu_\ell
=
r\exp\left(\frac{(2\ell+1)\pi i}{d}\right),
\qquad \ell=0,\ldots,d-1.
\]
\end{proof}

Since $R=I+W$, the eigenvalues of $R$ are
$\lambda_\ell=1+\nu_\ell$.
For $0\leq j<M$, recall that
\[
\theta_j=\frac{(2j+1)\pi}{d},
\qquad
\lambda_j=1+re^{i\theta_j}
=\rho_j e^{i\phi_j}.
\]
Since $0<\theta_j<\pi$, we have
$\operatorname{Im}\lambda_j=r\sin\theta_j>0$,
so $\phi_j\in(0,\pi)$. Moreover,
\[
\rho_j^2
=
|1+re^{i\theta_j}|^2
=
1+r^2+2r\cos\theta_j.
\]
Thus $\rho_j>1$ if and only if
$r+2\cos\theta_j>0$.
The remaining nonreal eigenvalues occur in
conjugate pairs, since
$\nu_{d-1-j}=\overline{\nu_j}$.
When $d$ is odd, the additional real eigenvalue
of $W$ is $\nu_M=-r$, yielding the real
eigenvalue $1-r$ of $R$.

\subsection{Diagonal similarity to a canonical wrapped shift}
\label{app:diagonal similarity to a canonical wrapped shift}
The spectrum above can also be obtained through an explicit diagonal similarity. This similarity shows that the spectrum of a wrapped weighted backward shift depends on the weights only through their full-cycle product.

Define the canonical wrapped weighted backward shift
$
C:\mathbb R^d\to\mathbb R^d
$
by
$
Ce_0=-r^de_{d-1}
$
and
$
Ce_j=e_{j-1},
\
j=1,\ldots,d-1.
$
Thus $C$ is the wrapped weighted backward shift whose weights are
$
(-r^d,1,\ldots,1).
$
Define nonzero scalars
$
\gamma_0:=1
$
and
\[
\gamma_j
:=
\frac{1}{w_1w_2\cdots w_j},
\qquad
j=1,\ldots,d-1.
\]
Let
$
D:\mathbb R^d\to\mathbb R^d
$
be the invertible diagonal operator defined by
$
De_j=\gamma_je_j.
$
 
\begin{lemma}[Explicit diagonal similarity]
\label{lem:diagonal similarity canonical shift}
One has
$
D^{-1}WD=C.
$
Equivalently,
$
W=DCD^{-1}.
$
\end{lemma}

\begin{proof}
For $1\leq j<d$, the definition of $\gamma_j$
gives
\[
D^{-1}WD e_j
=
\frac{\gamma_jw_j}{\gamma_{j-1}}e_{j-1}
=
e_{j-1}
=
Ce_j.
\]
For $j=0$, using
$\Omega=\prod_{j=0}^{d-1}w_j=-r^d$,
\[
D^{-1}WD e_0
=
\frac{w_0}{\gamma_{d-1}}e_{d-1}
=
\Omega e_{d-1}
=
Ce_0.
\]
Thus $D^{-1}WD=C$.
\end{proof}

The canonical shift $C$ has an explicit eigenvector associated with every root $\nu$ of
$
\nu^d=-r^d.
$
Indeed, define
$
v(\nu)
:=
(1,
\nu,
\nu^2,
\dots,
\nu^{d-1})^{\top}
\in\mathbb C^d.
$
For every root $\nu$ of $\nu^d=-r^d$,
\[
Cv(\nu)
=
(\nu,\nu^2,\ldots,\nu^d)^\top
=
\nu v(\nu).
\]
Since $W=DCD^{-1}$, it follows that
$z(\nu):=Dv(\nu)$ satisfies
\[
Wz(\nu)=WDv(\nu)=DCv(\nu)=\nu z(\nu).
\]


\subsection{Real rotation--dilation decomposition}
\label{app:Real rotation--dilation decomposition}
We next pass from the complex eigenvalue description to a real block decomposition.
For every
$
j=0,\ldots,M-1,
$
let
$
z_j:=z(\nu_j)\in\mathbb C^d
$
be a complex eigenvector satisfying
$
Wz_j=\nu_jz_j.
$
Write
$
z_j=a_j+ib_j,
\
a_j,b_j\in\mathbb R^d.
$
The vectors $a_j$ and $b_j$ are linearly independent over $\mathbb R$. Indeed, suppose that they were linearly dependent. Then $z_j$ would be a nonzero complex scalar multiple of a real vector $a\in\mathbb R^d$. Thus
$
z_j=ca
$
for some
$
c\in\mathbb C\setminus\{0\}.
$
Since
$
Wz_j=\nu_jz_j,
$
we would obtain
$
cWa
=
\nu_jca,
$
and therefore
$
Wa=\nu_ja.
$
The left-hand side is real, whereas $\nu_j$ is nonreal and $a\neq0$. This is impossible. Hence
$
\dim_{\mathbb R}\operatorname{span}_{\mathbb R}\{a_j,b_j\}=2.
$
Define the real plane
$
E_j
:=
\operatorname{span}_{\mathbb R}\{a_j,b_j\}.
$
Because
$
Rz_j=\lambda_jz_j
$
and
$
\lambda_j
=
\rho_j
\left(
\cos\phi_j+i\sin\phi_j
\right),
$
we have
$
R(a_j+ib_j)
=
\rho_j
\left(
\cos\phi_j+i\sin\phi_j
\right)
(a_j+ib_j).
$
Expanding the right-hand side gives
\[
R(a_j+ib_j)
=
\rho_j
\left[
\cos\phi_j\,a_j
-
\sin\phi_j\,b_j
\right]
+
i\rho_j
\left[
\sin\phi_j\,a_j
+
\cos\phi_j\,b_j
\right].
\]
Comparing real and imaginary parts, we obtain
$
Ra_j
=
\rho_j
\left(
\cos\phi_j\,a_j
-
\sin\phi_j\,b_j
\right)
$
and
$
Rb_j
=
\rho_j
\left(
\sin\phi_j\,a_j
+
\cos\phi_j\,b_j
\right).
$

If $d$ is odd, let
$
z_\ast:=z(-r)\in\mathbb R^d.
$
Then
$
Wz_\ast=-rz_\ast
$
and
$
Rz_\ast=(1-r)z_\ast.
$
Define
$
E_\ast:=\operatorname{span}_{\mathbb R}\{z_\ast\}.
$

Because the eigenvalues of $W$ are all distinct, the corresponding complex eigenvectors form a basis of $\mathbb C^d$. Grouping conjugate pairs and, when $d$ is odd, the remaining real eigenvector, we obtain the real direct-sum decomposition
\[
\mathbb R^d
=
\bigoplus_{j=0}^{M-1}E_j
\oplus
E_\ast^\delta.
\]
Here the last summand is omitted when $\delta=0$.
Define the real model space
\[
\widehat Z
:=
\left(
\bigoplus_{j=0}^{M-1}\mathbb R^2
\right)
\oplus
\mathbb R^\delta.
\]
For
\[
\xi
=
\left(
(x_0,y_0),\ldots,(x_{M-1},y_{M-1}),t
\right)
\in\widehat Z,
\]
where the coordinate $t$ is omitted when $\delta=0$, define
\begin{equation}
\label{eq:expression U}
\mathcal U\xi
:=
\sum_{j=0}^{M-1}
\left(
x_ja_j+y_jb_j
\right)
+
\delta\,t z_\ast.
\end{equation}
The preceding direct-sum decomposition shows that
$
\mathcal U:\widehat Z\to\mathbb R^d
$
is a real linear isomorphism.

\begin{proposition}[Real rotation--dilation similarity]
\label{prop:real rotational similarity R}
With the notation above,
\[
\mathcal U^{-1}R^n\mathcal U
=
\left(
\bigoplus_{j=0}^{M-1}
\rho_j^n\operatorname{Rot}(-n\phi_j)
\right)
\oplus
[1-r]^{n\delta}.
\]
\end{proposition}

\begin{proof}
On the $j$-th copy of $\mathbb R^2$, $\mathcal U_j(x, y)=x a_j+y b_j$. Then, $R\mathcal U_j(x, y)=xR a_j+yR b_j=\rho_j
(
x\cos\phi_j+y \sin \phi_j)\,a_j
+
\rho_j (-x\sin\phi_j+y \cos \phi_j)\,b_j.
$ 
Applying $\mathcal U^{-1}_j$ returns the coordinates
\[
\mathcal U^{-1}_j R\mathcal U_j  (x, y)= \rho_j
\begin{pmatrix}
\cos\phi_j&\sin\phi_j\\
-\sin\phi_j&\cos\phi_j
\end{pmatrix}
\begin{pmatrix}
x\\
y
\end{pmatrix}
=
\rho_j\operatorname{Rot}(-\phi_j) (x, y).
\]
Thus, on the $j$-th two-dimensional summand,
$
\mathcal U_j^{-1}R\mathcal U_j
=
\rho_j\operatorname{Rot}(-\phi_j).
$
If $d$ is odd, $\mathcal U$ maps the canonical basis of the final one-dimensional summand to $z_\ast$, and
$
Rz_\ast=(1-r)z_\ast.
$
Hence the final one-dimensional block is
$
[1-r].
$
Combining all invariant summands yields $\mathcal U^{-1}R\mathcal U
=
\left(
\bigoplus_{j=0}^{M-1}
\rho_j\operatorname{Rot}(-\phi_j)
\right)
\oplus
[1-r]^{\delta}$.

Moreover, $(\mathcal U^{-1}R\mathcal U)^2=\mathcal U^{-1}R\mathcal U \mathcal U^{-1}R\mathcal U=\mathcal U^{-1}R^2\mathcal U$. So, for every $n\ge 1, (\mathcal U^{-1}R\mathcal U)^n=\mathcal U^{-1}R^n\mathcal U$. Therefore, the corresponding diagram commutes:
\[
\begin{array}{ccc}
\mathbb R^d
&
\xrightarrow{\ R^n\ }
&
\mathbb R^d
\\
\mathcal U^{-1}\downarrow
&&
\downarrow\mathcal U^{-1}
\\
\widehat Z
&
\xrightarrow{
\left(
\bigoplus_{j=0}^{M-1}
\rho_j^n\operatorname{Rot}(-n\phi_j)
\right)
\oplus
[1-r]^{n\delta}
}
&
\widehat Z.
\end{array}
\]

To justify the lower horizontal arrow, observe that
$
\operatorname{Rot}(\phi)\operatorname{Rot}(\psi)
=
\operatorname{Rot}(\phi+\psi).
$
Hence, by induction,
$
\operatorname{Rot}(\phi)^n
=
\operatorname{Rot}(n\phi).
$
Therefore
$
\left(
\rho_j\operatorname{Rot}(-\phi_j)
\right)^n
=
\rho_j^n\operatorname{Rot}(-n\phi_j).
$
\end{proof}

The purpose of $\mathcal U$ is to provide a fixed real coordinate system in which the conjugate spectral modes become explicit rotation--dilation planes. We will use these coordinates to define the source, control, and output channels on the original space.

\subsection{Source and control blocks.}
\label{app:explicit construction channels for R}
Choose a set of conjugate-plane indices
$
J
\subseteq
\{0,\ldots,M-1\}.
$
Write
$
J=\{j_1,\ldots,j_p\},
$
where
$
p:=|J|
\leq
\left\lfloor\frac d2\right\rfloor.
$
Let
$
Y:=\mathbb R^p
$
with canonical basis
$
f_1,\ldots,f_p.
$
For every selected plane $E_{j_\ell}$, we use: the first coordinate of the plane as the embedded source coordinate;
 the second coordinate of the plane as the control coordinate;
and the first coordinate of the plane as the output readout.

We first define the channels on the model space $\widehat Z$. Let
$
\widehat J_{\mathrm{out}}:Y\to\widehat Z
$
be defined by
\begin{equation}
\label{eq:expression j hat out}
\widehat J_{\mathrm{out}}f_\ell
=
\left(
0,\ldots,0,
\underset{j_\ell\text{-th plane}}{(1,0)},
0,\ldots,0
\right).
\end{equation}
Likewise, define
$
\widehat J_{\mathrm{ctl}}:Y\to\widehat Z
$
by
\[
\widehat J_{\mathrm{ctl}}f_\ell
=
\left(
0,\ldots,0,
\underset{j_\ell\text{-th plane}}{(0,1)},
0,\ldots,0
\right).
\]

Define
$
\widehat P:\widehat Z\to Y
$
by
\[
\widehat P
\left(
(x_0,y_0),\ldots,(x_{M-1},y_{M-1}),t
\right)
=
\left(
x_{j_1},
\ldots,
x_{j_p}
\right),
\]
where the final coordinate $t$ is omitted when $d$ is even.

Thus $\widehat P$ reads the first coordinate of every selected conjugate plane and ignores:
 the second coordinate of every selected plane;
 both coordinates of every unselected plane;
and the remaining real eigendirection when $d$ is odd.

The model-space channels satisfy
$
\widehat P\widehat J_{\mathrm{out}}
=
I_Y
$
and
$
\widehat P\widehat J_{\mathrm{ctl}}
=
0.
$
Indeed, for every canonical basis vector $f_\ell$,
$
\widehat P\widehat J_{\mathrm{out}}f_\ell
=
f_\ell,
$
whereas
$
\widehat P\widehat J_{\mathrm{ctl}}f_\ell
=
0.
$

We now transport these channels to the original state space
$
Z=\mathbb R^d.
$
Define
$
J_{\mathrm{out}}
:=
\mathcal U\widehat J_{\mathrm{out}},
\
J_{\mathrm{ctl}}
:=
\mathcal U\widehat J_{\mathrm{ctl}},
$
and
$
P
:=
\widehat P\mathcal U^{-1}.
$
Then
$
J_{\mathrm{out}},J_{\mathrm{ctl}}:Y\to Z
$
are injective, while
$
P:Z\to Y
$
is surjective. Moreover,
\[
PJ_{\mathrm{out}}
=
\widehat P\mathcal U^{-1}
\mathcal U\widehat J_{\mathrm{out}}
=
\widehat P\widehat J_{\mathrm{out}}\\
=
I_Y,
\]
and
\[
PJ_{\mathrm{ctl}}
=
\widehat P\mathcal U^{-1}
\mathcal U\widehat J_{\mathrm{ctl}}
=
\widehat P\widehat J_{\mathrm{ctl}}
=
0.
\]
Thus these are legitimate output and control channels in the notation of the Steering Criterion.
For every $n\geq0$, define
$
A_n:=PR^nJ_{\mathrm{out}}
$
and
$
G_n:=PR^nJ_{\mathrm{ctl}}.
$

\begin{lemma}[Exact source and control blocks]
\label{lem:rotational source control blocks}
For every $n\geq0$,
$$
A_n
=
\operatorname{diag}_{1\leq\ell\leq p}
\left(
\rho_{j_\ell}^n
\cos(n\phi_{j_\ell})
\right)
\quad
\text{ and }
\quad
G_n
=
\operatorname{diag}_{1\leq\ell\leq p}
\left(
\rho_{j_\ell}^n
\sin(n\phi_{j_\ell})
\right).
$$
\end{lemma}

\begin{proof}
Fix
$
\ell\in\{1,\ldots,p\}.
$
By construction,
$
\widehat J_{\mathrm{out}}f_\ell
$
is the first canonical coordinate vector
$
(1,0)
$
in the plane indexed by $j_\ell$, and
$
J_{\mathrm{out}}f_\ell
=
\mathcal U\widehat J_{\mathrm{out}}f_\ell.
$
By Proposition \ref{prop:real rotational similarity R},
we have $
\mathcal U^{-1}R^n\mathcal U
=
\left(
\bigoplus_{j=0}^{M-1}
\rho_j^n\operatorname{Rot}(-n\phi_j)
\right)
\oplus
[1-r]^{n\delta},
$
so we obtain
$
A_nf_\ell
=
PR^nJ_{\mathrm{out}}f_\ell
=
\widehat P\mathcal U^{-1}
R^n
\mathcal U
\widehat J_{\mathrm{out}}f_\ell.
$
On the selected plane $j_\ell$,
\[
\rho_{j_\ell}^n
\operatorname{Rot}(-n\phi_{j_\ell})
\begin{pmatrix}
1\\
0
\end{pmatrix}
=
\rho_{j_\ell}^n
\begin{pmatrix}
\cos(n\phi_{j_\ell})\\
-\sin(n\phi_{j_\ell})
\end{pmatrix}.
\]
The map $\widehat P$ reads the first coordinate of this plane. Therefore
$
A_nf_\ell
=
\rho_{j_\ell}^n
\cos(n\phi_{j_\ell})
f_\ell.
$
Since this holds for every $\ell$,
$
A_n
=
\operatorname{diag}_{1\leq\ell\leq p}
\left(
\rho_{j_\ell}^n
\cos(n\phi_{j_\ell})
\right).
$

Similarly,
$
\widehat J_{\mathrm{ctl}}f_\ell
$
is the second canonical coordinate vector
$
(0,1)
$
in the plane indexed by $j_\ell$. On that plane,
\[
\rho_{j_\ell}^n
\operatorname{Rot}(-n\phi_{j_\ell})
\begin{pmatrix}
0\\
1
\end{pmatrix}
=
\rho_{j_\ell}^n
\begin{pmatrix}
\sin(n\phi_{j_\ell})\\
\cos(n\phi_{j_\ell})
\end{pmatrix}.
\]
Reading the first coordinate gives
$
G_nf_\ell
=
\rho_{j_\ell}^n
\sin(n\phi_{j_\ell})
f_\ell.
$
Thus,
$
G_n
=
\operatorname{diag}_{1\leq\ell\leq p}
\left(
\rho_{j_\ell}^n
\sin(n\phi_{j_\ell})
\right).
$
\end{proof}

The operator $G_n$ is invertible if and only if
$
\sin(n\phi_{j_\ell})\neq0
\
\text{for every }
\ell=1,\ldots,p.
$
Whenever this holds,
$
G_n^{-1}
=
\operatorname{diag}_{1\leq\ell\leq p}
\left(
\rho_{j_\ell}^{-n}
\csc(n\phi_{j_\ell})
\right).
$
Furthermore,
\[
\begin{aligned}
G_n^{-1}A_n
&=
\operatorname{diag}_{1\leq\ell\leq p}
\left(
\rho_{j_\ell}^{-n}
\csc(n\phi_{j_\ell})
\right)
\operatorname{diag}_{1\leq\ell\leq p}
\left(
\rho_{j_\ell}^n
\cos(n\phi_{j_\ell})
\right)
=
\operatorname{diag}_{1\leq\ell\leq p}
\left(
\csc(n\phi_{j_\ell})
\cos(n\phi_{j_\ell})
\right)\\
&=
\operatorname{diag}_{1\leq\ell\leq p}
\left(
\cot(n\phi_{j_\ell})
\right).
\end{aligned}
\]
Thus the source-correction maps and target-control maps determined by the exact block characterization are
\begin{equation}
\label{eq:inverse diag control formulas for R}
Q_n
=
-G_n^{-1}A_n
=
-\operatorname{diag}_{1\leq\ell\leq p}
\left(
\cot(n\phi_{j_\ell})
\right),
\quad
S_n
=
G_n^{-1}
=
\operatorname{diag}_{1\leq\ell\leq p}
\left(
\rho_{j_\ell}^{-n}
\csc(n\phi_{j_\ell})
\right).
\end{equation}

We equip
$
Y=\mathbb R^p
$
with its Euclidean norm. Hence the operator norm of a diagonal operator
$
\operatorname{diag}(c_1,\ldots,c_p)
$
is
$
\max_{1\leq\ell\leq p}|c_\ell|.
$
Therefore
$
\|G_n^{-1}A_n\|
=
\max_{1\leq\ell\leq p}
|\cot(n\phi_{j_\ell})|
$
and
$
\|G_n^{-1}\|
=
\max_{1\leq\ell\leq p}
\left[
\rho_{j_\ell}^{-n}
|\csc(n\phi_{j_\ell})|
\right].
$
The Steering Criterion for these channels is therefore equivalent to finding an increasing sequence $(n_k)$ such that simultaneously
$
\cot(n_k\phi_{j_\ell})\longrightarrow0
$
and
$
\rho_{j_\ell}^{-n_k}
|\csc(n_k\phi_{j_\ell})|
\longrightarrow0
$
for every
$
\ell=1,\ldots,p.
$
The first family of conditions is purely rotational. The second combines rotational recurrence with radial expansion.

\subsection{The torus condition}
\label{app:torus condition}
Use \(\alpha_J\), \(\beta\), and \(H_{\alpha_J}\) as defined
in Section~\ref{sec:steering for R}.

\begin{lemma}[Phase convergence and cotangent collapse]
\label{lem:phase convergence cotangent}
Let $(t_k)$ be a sequence of real numbers. Then
$
\cot(t_k)\longrightarrow0
$
if and only if
$
\frac{t_k}{\pi}
\longrightarrow
\frac12
\
\text{in }\mathbb T.
$
In either case,
$
|\sin t_k|\longrightarrow1
$
and
$
|\csc t_k|\longrightarrow1.
$
\end{lemma}

\begin{proof}
The equivalence follows from the continuity
and $\pi$-periodicity of the cotangent,
whose zeros are precisely $\pi/2$ modulo
$\pi$. The remaining limits follow from
$1+\cot^2t=\csc^2t$.
\end{proof}


The following is the standard minimal-rotation property of the
compact cyclic group $H_{\alpha_J}$; see
\citep[Theorem~10.13]{eisner2015}.
\begin{lemma}[Existence of an increasing phase-recurrence sequence]
\label{lem:increasing torus recurrence sequence}
There exists a strictly increasing sequence
$
n_k\longrightarrow\infty
$
such that
$
n_k\alpha_J\longrightarrow\beta
\
\text{in }\mathbb T^p
$
if and only if
$
\beta\in H_{\alpha_J}.
$
\end{lemma}



\subsection{Detailed proof of Theorem \ref{thm:exact spectral toral classification R} (Exact spectral–toral classification for $R$)}
\label{app:detailed proof of exact spectral toral classification for R}
The equivalence between {\rm(ii)} and {\rm(iii)} follows from Kronecker's closure theorem \citep{Cassels1957} and the parity-resonance calculation above.

\textit{(ii) implies (i).}
Assume that $\rho_{j_\ell}>1$ for every
$\ell=1,\ldots,p$ and $\beta\in H_{\alpha_J}$.
By Lemma~\ref{lem:increasing torus recurrence sequence}, there exists a strictly increasing
sequence $n_k\to\infty$ such that
$n_k\alpha_J\to\beta$ in $\mathbb T^p$.
Lemma~\ref{lem:phase convergence cotangent} then gives
\[
\max_{1\leq\ell\leq p}
|\cot(n_k\phi_{j_\ell})|
\longrightarrow0,
\qquad
|\csc(n_k\phi_{j_\ell})|
\longrightarrow1
\quad\text{for every }\ell.
\]
In particular, $G_{n_k}$ is eventually
invertible. By Lemma~\ref{lem:rotational source control blocks},
\[
\|G_{n_k}^{-1}A_{n_k}\|
=
\max_{1\leq\ell\leq p}
|\cot(n_k\phi_{j_\ell})|
\longrightarrow0.
\]
Since all selected radial factors exceed one,
\[
\|G_{n_k}^{-1}\|
=
\max_{1\leq\ell\leq p}
\rho_{j_\ell}^{-n_k}
|\csc(n_k\phi_{j_\ell})|
\longrightarrow0.
\]
Theorem~\ref{thm:block charact} therefore gives (i).

\medskip
\noindent
\textit{(i) implies (ii).}
Assume that the Steering Criterion holds
along a strictly increasing sequence $n_k\to\infty$.
By Theorem~\ref{thm:block charact},
$
\|G_{n_k}^{-1}\|\longrightarrow0,
$ and $
\|G_{n_k}^{-1}A_{n_k}\|
\longrightarrow0.
$
Lemma~\ref{lem:rotational source control blocks} shows that
$\cot(n_k\phi_{j_\ell})\to0$ for every
$\ell$. By Lemma~\ref{lem:phase convergence cotangent},
$n_k\alpha_J\to\beta$ in $\mathbb T^p$,
so $\beta\in H_{\alpha_J}$, and
$
|\csc(n_k\phi_{j_\ell})|
\longrightarrow1
\quad\text{for every }\ell.
$
Since
\[
\rho_{j_\ell}^{-n_k}
|\csc(n_k\phi_{j_\ell})|
\leq\|G_{n_k}^{-1}\|
\longrightarrow0,
\]
we must have $\rho_{j_\ell}>1$ for every
selected mode. Thus (ii) holds.

\emph{{\rm(ii)} implies {\rm(iv)}.}
Assume that
$
\rho_{j_\ell}>1
\
\forall\ell=1,\ldots,p,
$
and that
$
\beta\in H_{\alpha_J}.
$
Fix $\varepsilon>0$ and define the relative neighborhood
\[
U_\varepsilon
:=
\left\{
x\in H_{\alpha_J}:
\|x-\beta\|_{\mathbb T^p}<\varepsilon
\right\}.
\]
Since
$
\beta\in H_{\alpha_J},
$
the set $U_\varepsilon$ is a nonempty relatively open subset of
$H_{\alpha_J}$.
Consider the translation
\[
\tau_{\alpha_J}:H_{\alpha_J}\to H_{\alpha_J},
\qquad
\tau_{\alpha_J}(x):=x+\alpha_J.
\]

Since $H_{\alpha_J}$ is the compact cyclic
group generated by $\alpha_J$, the translation
$\tau_{\alpha_J}$ is minimal
\citep[Theorem~10.13]{eisner2015}.
By the classical uniform-recurrence theorem
\citep{furstenberg1981},
the return-time set
\[
N_\varepsilon
=
\{n\geq 1:
\tau_{\alpha_J}^n(0)\in U_\varepsilon\}
\]
is syndetic. Thus (iv) holds.

\emph{{\rm(iv)} implies {\rm(ii)}.}
Assume condition {\rm(iv)}. The radial condition
$
\rho_{j_\ell}>1
\
\forall\ell=1,\ldots,p
$
is already part of {\rm(iv)}. It remains to prove that
$
\beta\in H_{\alpha_J}.
$
For every integer $k\geq1$, set
$
\varepsilon_k:=\frac1k.
$
Suppose we have already chosen $n_{k-1}\in \mathcal N_{\varepsilon_{k-1}}$. Since $\mathcal N_{\varepsilon_k}$ is syndetic, hence infinite, we can choose
$
n_k\in\mathcal N_{\varepsilon_k}
$ such that $n_{k-1}<n_k$.
Then we have an increasing sequence $(n_k)$ such that
$
\|n_k\alpha_J-\beta\|_{\mathbb T^p}
<
\frac1k,
$ for all $k$.
Consequently, 
$
n_k\alpha_J
\longrightarrow
\beta
\
\text{in }\mathbb T^p$. By Lemma \ref{lem:increasing torus recurrence sequence}, 
$
\beta
\in
H_{\alpha_J}.
$
Thus condition {\rm(ii)} holds.
\hfill$\blacksquare$

\subsection{Smallest residual radial factor for expansivity}
\label{app:smallest residual radial factor}

\begin{lemma}[Smallest residual radial factor]
\label{lem:smallest residual radial factor}
The sequence
$
(\rho_j)_{j=0}^{M-1}
$
is strictly decreasing. Consequently,
$
\min_{0\leq j\leq M-1}\rho_j
=
\rho_{M-1}.
$
Moreover,
$
\theta_{M-1}
=
\pi-\frac{\kappa_d\pi}{d},
$
and hence
$
\rho_{M-1}^2
=
1+r^2
-
2r\cos\left(
\frac{\kappa_d\pi}{d}
\right).
$
Therefore all conjugate residual modes are expansive if and only if
$
r
>
2\cos\left(
\frac{\kappa_d\pi}{d}
\right).
$
\end{lemma}

\begin{proof}
Since
$
0<\theta_0<\cdots<\theta_{M-1}<\pi
$
and cosine is strictly decreasing on $(0,\pi)$,
the identity
$
\rho_j^2=1+r^2+2r\cos\theta_j
$
shows that $(\rho_j)_{j=0}^{M-1}$
is strictly decreasing. Hence
$
\min_{0\leq j<M}\rho_j=\rho_{M-1}.
$
Since $d=2M+\delta$, with
$\delta\in\{0,1\}$, and
$\kappa_d=1+\delta$, we have
\[
\theta_{M-1}
=
\frac{(2M-1)\pi}{d}
=
\pi-\frac{\kappa_d\pi}{d}.
\]
Consequently,
\[
\rho_{M-1}^2
=
1+r^2-2r\cos\left(\frac{\kappa_d\pi}{d}\right).
\]
All conjugate residual modes are expansive
if and only if $\rho_{M-1}>1$, equivalently,
since $r>0$, to
$
r>2\cos\left(\frac{\kappa_d\pi}{d}\right).
$
\end{proof}

\section{Obstruction to Asymptotic Full-State Steering}
\label{app:full-state-obstruction}

\begin{proposition}[Obstruction to asymptotic full-state steering]
\label{prop:full-state-obstruction}
Let $T\in\mathcal L(Z)$ be in either of the following settings:

\noindent
\textup{(i)} $T=W$ is a wrapped weighted backward shift satisfying $|\Omega|>1$.

\noindent
\textup{(ii)} $T=R=I+W$, with selected conjugate planes $J=\{j_1,\ldots,j_p\}$ satisfying $\rho_{j_\ell}>1$ for every $\ell$, and with the corresponding channels $J_{\mathrm{out}}$ and $J_{\mathrm{ctl}}$.

Let $(n_k)$ be any strictly increasing sequence and let $Q_k,S_k:Y\to Y$ satisfy $\|Q_k\|\to0$ and $\|S_k\|\to0$. For $u,v\in Y$, define
$x_k(u,v):=J_{\mathrm{out}}u+J_{\mathrm{ctl}}(Q_ku+S_kv)$.

Then, for every $u\neq0$ and every $v\in Y$,
\[
\|T^{n_k}x_k(u,v)-J_{\mathrm{out}}v\|
\longrightarrow\infty.
\]
In particular, neither steering mechanism can be asymptotically full-state accurate.
\end{proposition}

\begin{proof}
Fix $u\neq0$ and $v\in Y$. Since $\|Q_k\|,\|S_k\|\to0$, we have
$x_k(u,v)\to J_{\mathrm{out}}u$.

The argument uses the same expanding-subspace principle in both cases. We identify a $T$-invariant subspace $E\subseteq Z$ and a projection $\Pi:Z\to E$ commuting with $T$ such that
$\Pi J_{\mathrm{out}}=J_{\mathrm{out}}$ and
$\|(T|_E)^{-n}\|\to0$.

In case \textup{(i)}, take $E=Z$ and $\Pi=I_Z$. Writing $n=ad+r$, with $0\leq r<d$, the cycle identity $W^d=\Omega I_Z$ gives
$W^{-n}=\Omega^{-a}W^{-r}$.
Consequently,
$\|W^{-n}\|\leq C|\Omega|^{-a}\to0$,
where $C:=\max_{0\leq r<d}\|W^{-r}\|$.

In case \textup{(ii)}, let
$E:=\bigoplus_{\ell=1}^p E_{j_\ell}$ defined as in Appendix \ref{app:Real rotation--dilation decomposition}
be the direct sum of the selected invariant conjugate planes. We established that the restriction $R|_E$ is similar to
$\bigoplus_{\ell=1}^p\rho_{j_\ell}\operatorname{Rot}(-\phi_{j_\ell})$.
Hence, setting
$\rho_{\min}:=\min_{1\leq\ell\leq p}\rho_{j_\ell}>1$,
there exists $C>0$ such that
$\|(R|_E)^{-n}\|\leq C\rho_{\min}^{-n}\to0$.
The projection $\Pi$ onto $E$ along the remaining invariant spectral subspaces commutes with $R$. Moreover, the construction of the channels gives
$J_{\mathrm{out}}Y\subseteq E$ (by definition of $J_{\rm out}$ and Equations (\ref{eq:expression U}) and (\ref{eq:expression j hat out})), so
$\Pi J_{\mathrm{out}}=J_{\mathrm{out}}$.

Thus, in either case,
$\Pi x_k(u,v)\to J_{\mathrm{out}}u\neq0$,
where the last inequality follows from $PJ_{\mathrm{out}}=I_Y$.
Since $\Pi$ commutes with $T$, we obtain
\[
\begin{aligned}
\|\Pi\|\,\|T^{n_k}x_k(u,v)\|
&\geq
\|T^{n_k}\Pi x_k(u,v)\|\\
&\geq
\frac{\|\Pi x_k(u,v)\|}
{\|(T|_E)^{-n_k}\|}
\longrightarrow\infty.
\end{aligned}
\]
Therefore $\|T^{n_k}x_k(u,v)\|\to\infty$. Finally, the reverse triangle inequality yields
\[
\|T^{n_k}x_k(u,v)-J_{\mathrm{out}}v\|
\geq
\|T^{n_k}x_k(u,v)\|-\|J_{\mathrm{out}}v\|
\longrightarrow\infty.
\]
This proves the obstruction in both cases.
\end{proof}

\section{Proof of Equation (\ref{eq:exact preservation of the same realization across depths})}
\label{app:Contragredient Coefficient Transport}

\begin{proposition}
\label{prop:contragredient coefficient transport}
For every $c\in Y_{\mathrm{fr}}$ and every $s\geq0$,
\[
\mathscr D_s
\left(
\mathcal C_{\mathrm{scaff}}^sc
\right)
=
\mathscr D_0(c).
\]
Moreover,
$\mathcal C_{\mathrm{scaff}}^{d_{\mathrm{scaff}}}
=\Omega_{\mathrm{scaff}}^{-1}I_{Y_{\mathrm{fr}}}$.
Consequently, if
$s=ad_{\mathrm{scaff}}+t$,
where $0\leq t<d_{\mathrm{scaff}}$, then
$\mathcal C_{\mathrm{scaff}}^s
=\Omega_{\mathrm{scaff}}^{-a}\mathcal C_{\mathrm{scaff}}^t$.
\end{proposition}

\begin{proof}
We first establish the one-step identity
$\mathscr D_1(\mathcal C_{\mathrm{scaff}}c)
=\mathscr D_0(c)$.
By definition,
\[
\begin{aligned}
\mathscr D_1(\mathcal C_{\mathrm{scaff}}c)
=
\sum_{k=1}^K
\sum_{j=1}^{|D_k|}
\sum_{m=0}^{d_{\mathrm{scaff}}-2}
c_{k,m+1,j}
J_kW_{\mathrm{scaff}}^{m+1}x_{k,j}
+
\sum_{k=1}^K
\sum_{j=1}^{|D_k|}
\Omega_{\mathrm{scaff}}^{-1}c_{k,0,j}
J_kW_{\mathrm{scaff}}^{d_{\mathrm{scaff}}}x_{k,j}.
\end{aligned}
\]
Since
$W_{\mathrm{scaff}}^{d_{\mathrm{scaff}}}
=\Omega_{\mathrm{scaff}}I_X$,
the second sum equals
$\sum_{k,j}c_{k,0,j}J_kx_{k,j}$.
Relabeling the first sum therefore gives
\[
\mathscr D_1(\mathcal C_{\mathrm{scaff}}c)
=
\sum_{k=1}^K
\sum_{m=0}^{d_{\mathrm{scaff}}-1}
\sum_{j=1}^{|D_k|}
c_{k,m,j}J_kW_{\mathrm{scaff}}^mx_{k,j}
=
\mathscr D_0(c).
\]

The same calculation, with $W_{\mathrm{scaff}}^{m+s}$ in place of $W_{\mathrm{scaff}}^m$, yields
$\mathscr D_{s+1}(\mathcal C_{\mathrm{scaff}}c)
=\mathscr D_s(c)$
for every $s\geq0$.
Iterating this identity gives
$\mathscr D_s(\mathcal C_{\mathrm{scaff}}^sc)
=\mathscr D_0(c)$.

The operator $\mathcal C_{\mathrm{scaff}}$ shifts each depth coefficient one position backward in its indexed coefficient cycle, multiplying it by
$\Omega_{\mathrm{scaff}}^{-1}$
when it crosses the wrap. After exactly $d_{\mathrm{scaff}}$ steps, each coefficient has crossed the wrap once. Hence
$\mathcal C_{\mathrm{scaff}}^{d_{\mathrm{scaff}}}
=\Omega_{\mathrm{scaff}}^{-1}I_{Y_{\mathrm{fr}}}$.
Writing $s=ad_{\mathrm{scaff}}+t$ gives the final identity.
\end{proof}

\end{document}